\documentclass[preprint]{elsarticle}

\usepackage{hyperref}

\journal{Journal of \LaTeX\ Templates}

\usepackage{amsmath}
\usepackage{amsfonts}
\usepackage{amssymb}
\usepackage{graphicx,xcolor}

\usepackage[english]{babel}
\usepackage[utf8x]{inputenc}
\usepackage[T1]{fontenc}

\def \R{\mbox{l\hspace{-.15em}R}}

\begin{document}

		\begin{frontmatter}
			\title{Padé-type High-Order Absorbing Boundary Condition for a Coupled Hydrodynamic Wave Model\\with Surface Tension Effect}
		
		\author{O. Wilk} \ead{olivier.wilk@cnam.fr}
		
		\address{Department of Mathematics\\
			Modélisation Mathématique et Numérique \\
			Conservatoire National des Arts et Métiers \\
			292, rue Saint Martin, Paris, 75003, France}
		
		\begin{abstract}
			This paper presents a specific way to develop High-Order Absorbing Boundary Conditions (HABC) of the Padé family on a Coupled Hydrodynamic Wave Model (CHWM) especially with surface tension effect (small scales in space). Inspired by the Neumann-Kelvin model, the CHWM is composed by a fluid-basin model to allow to use multiple objects below the surface coupled to a free surface model with a small added mass surface term. With the surface tension effect, we introduce new coefficients (similar as Higdon coefficients) on each HABC (for the surface model and the basin model) to ensure the continuity of the two HABC at the interface between the coupled models. So, we propose a useful specific compatibility condition, and a strong reduction of the Padé approximation in particular in the water case.
        \end{abstract}

		
		\end{frontmatter}
		
	
\section{Introduction}

In the context of numerical hydrodynamic wave simulations with multiple small objects just below the surface useful for example for fluid-structure interaction problems, we are interested in extending the possibilities to compute faster simulations for small or medium scales in space i.e. with the surface tension effect. We propose to take advantage of the performance of the high order absorbing boundary conditions to facilitate and control the use of the space below the surface. To do this, we choose the way to reduce the computational cost by using a truncated domain, the fluid or basin domain. Therefore, we need to develop a specific high order absorbing boundary condition just at the interface of the surface, which is the main purpose of this article. \\

Usually, the hydrodynamic problems are based on a model generated with an irrotational  incompressible fluid below the surface as the Neumann-Kelvin Model \cite{Stoker}. To explore such problems with the numerical simulations, these properties allow to reduce to a surface problem with Boundary Element Methods (BEM) \cite{BEM} \cite{Clement}	 \cite{VanDerStoep}. But if we want to study problems with objects just below the surface, the size and the number of these objects (as in \cite{melange-spheres-Valettea}, \cite{microcapsules-DeVuyst}) increase the numerical cost of the BEM simulations. It can be expensive. In this situation, the classical Finite Element Method (FEM) becomes interesting again, more competitive. With FEM, the fluid problem in a truncated domain (the basin domain) below the surface must be solved. But the associated incompressible model is not so useful to allow the wave outputs. Partial solutions use bigger domains or smaller study times. To solve this problem, we propose to forget the usual incompressibility of the fluid. We prefer to use the isentropic compressibility, the pressure, the density, the velocity evolve with the fluid celerity with the D'Alembert's equation, the wave equation \cite{Landau-Lifshitz}. This solution has been proposed \cite{Bissengue} and tested with the ABC of the order one. We use here the model without flow. \\

Also, we offer the opportunity to compute problems with small uniform added mass on the surface. We add an appropriate inertial term to the model. Now, our model can be perceived as a wave model (dispersive on the surface), the Coupled Hydrodynamic Wave Model (CHWM)  with the surface tension effect. This allows to explore the different families of Absorbing Boundary Conditions (ABCs) to build a specific ABC useful to truncate the domain. There are the Engquist-Majda ABCs \cite{Engquist-Majda-1977}, the Bayliss-Turkel ABCs \cite{Bayliss-Turkel-1980}, the High order ABCs (HABCs) associated to the approximation of the Dirichlet-to-Neumann map as the Padé-type HABCs (\cite{Engquist-Majda-1977}, \cite{Keller-Givoli-1989}, \cite{Collino-1993}), the HABCs of the family of the Higdon-HABCs (\cite{Higdon-1986}, \cite{Hagstrom-Hariharan-1998}, \cite{Givoli-Neta-2003}, \cite{Hagstrom-Warburton-2004}) and the Perfectly Matched Layers (PMLs) (\cite{Berenger-1994}, \cite{Turkel-Yefet-1998}). \\

\textbf{Remark:} The term "CHWM" \footnote{\textit{Most used abbreviations:}\\
	\begin{tabular}{l|l}
		~~ HWM: Hydrodynamic Wave Model & ABC: Absorbing Boundary Condition \\
		~~ CHWM: our Coupled HWM & HABC: High-Order ABC \\
		~~ STE: Surface Tension Effect & CP-HABC: Coupled Padé-type HABC
\end{tabular}}, in particular with the term "Coupled", is not perfectly judicious without the Surface Tension Effect (STE). The problem is local on the surface with no need for boundary condition. However, this term makes sense with the STE. The surface model needs a boundary condition, a specific compatible ABC.\\

With our CHWM without flow, we choose to apply the methodology associated to the Dirichlet-to-Neumann map with the Padé approximation in the case of the propagative waves. We can use reduced models either for the basin variable (the velocity potential) or the surface variable (the vertical displacement of the surface). So we get the Dirichlet to Neumann operator (with the square root applied to a part of each reduced models). This gives the  exact transparent boundary condition, a none local operator not useful for numerical simulations. We therefore use the classically Padé approximation to obtain a specific Padé-type HABC which is more interesting for numerical computations \cite{Engquist-Majda-1977}. To get efficient solutions in particular with the STE and the water case, we propose to apply a specific compatibility condition to ensure the continuity between the HABC at the boundaries of the surface problem and the HABC associated to the basin problem. This compatibility condition uses new coefficients in the spirit of the Higdon coefficients \cite{Higdon-1986}.\\

Following this introduction, we present in the next section the CHWM. In the section 3, we do some useful recalls associated the classical Padé-type HABC with the wave equation. In the section 4, we get the Padé-type HABC with the CHWM in the two cases, without and with the STE. In the first case, it's the easier model to get the Padé-type HABC, the surface model is local. The problem doesn't use boundary condition on the surface. The more interesting part (section 4.2) is associated to the CHWM with the STE, the surface model is not longer local. In this case, we must use two boundary conditions, one on the boundaries of basin domain and the other one on the boundaries of the surface domain. To get an efficient HABC, we add a compatibility condition. In the section 5, we propose a size reduction of the  number of the terms of the Padé approximation when we use an important value for the wave celerity of the fluid, for example with the water. The last section propose some numerical simulations without and with the STE for different values of the wave celerity of the fluid to evaluate the performance of our Padé-type HABC associated to our CHWM. We test also the size reduction of the number of the terms of the Padé approximation for the big values of the wave celerity of the fluid. And we finish with a case with a small object just below the surface.

\section{Hydrodynamic Wave Model}

\begin{figure}[h]
	\centering
	\includegraphics[width=8cm]{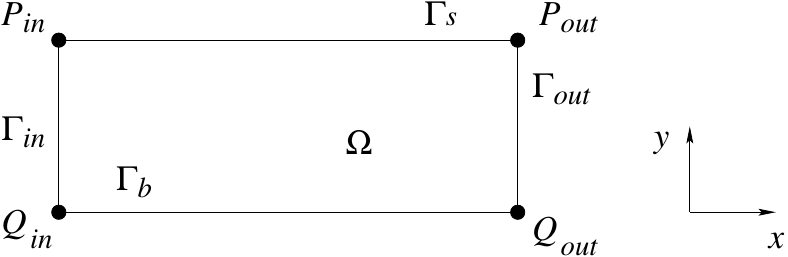}
	\caption{The domain and the boundaries.}
	\label{fig:domain}
\end{figure}
The HWM with or without STE is affiliated with the Neumann-Kelvin model here without flow. It's a coupled model with a surface model associated to the bernoulli equation and a potential fluid model (or basin model) below the surface $\Gamma_s$ in the domain $\Omega$ (Figure \ref{fig:domain}). The useful variables are $\eta$ the vertical displacement of the surface and $\varphi$ the velocity potential. A boundary condition at the surface connects the two models respecting the continuity of the normal velocity. So we get as first coupled model in the domain $\Omega$ (Figure \ref{fig:domain} and $\sigma, \rho, g$ respectively the surface tension coefficient, the mass density, the gravity, $\nu$ boundary normal, $t \in [0,T]$):
\begin{equation} 
	\left\lbrace
	\begin{array}{l}
		\displaystyle{ \sigma \partial_x^2 \eta  + \rho g \eta  + \rho \partial_t \varphi = f_s \mbox{ on } \Gamma_s,} \\
		\displaystyle{ \Delta \varphi = 0 \mbox{ in } \Omega, } \\
		\displaystyle{ \partial_{\nu} \varphi = \partial_t \eta \mbox{ on } \Gamma_s},\\
		\mbox{with boundary conditions on } \Gamma_{in} \cup \Gamma_{out} \cup  \Gamma_{b} \\
		\mbox{and initial conditions.}
	\end{array}
	\right.
\end{equation}
The surface model uses also a coupling term as well as the basin model with the Neumann condition on the boundary $\Gamma_s$ (continuity of the normal velocity). We wish to modify (to upgrade) the model with an added mass term in the surface model (for example with oil on water surface). And as explained in the introduction (paragraph 2) with the presence of several objects just below the surface so with FEM simulations, the usual incompressible fluid property is not so useful for wave outputs with a truncated domain. An isentropic compressible model \cite{Landau-Lifshitz}, with the wave equation, is more efficient for that. And eventually to test also faster problems, impact problems. So we replace the inner problem by the wave model. This coupled model was proposed in \cite{Bissengue} following \cite{Destuynder-Fabre-Wilk-2010} and \cite{Destuynder-Fabre-Orellana}. It is interesting, and not surprising, to note that the behavior of the CHWM is similar to the Neumann-Kelvin problem at the surface and also just below the surface especially with an important celerity coefficient $c_f$ as for the water. \\

\textbf{Remark:} For simple shape objects as cylinder, the presence of the object can be approximated with a Lamb’s submerged dipole \cite{Lamb} \cite{Moreira-peregrine} useful with BEM. But for more complex shapes, several objects and to study fluid-structure interaction, FEM becomes interesting again.\\

So we get our CHWM ($\varepsilon$ small positive coefficient associated the added mass term on the surface):
\begin{equation} \label{HWMeq}
\left\lbrace
\begin{array}{l}
	\displaystyle{ \varepsilon \rho \partial_t^2 \eta - \sigma \partial_x^2 \eta  + \rho g \eta  + \rho \partial_t \varphi = f_s \mbox{ on } \Gamma_s,} \\
	\displaystyle{ \partial_t^2 \varphi - c_f^2 \Delta \varphi = 0 \mbox{ in } \Omega, } \\
	\displaystyle{ \partial_{\nu} \varphi = \partial_t \eta \mbox{ on } \Gamma_s},\\
	\mbox{with boundary conditions on } \Gamma_{in} \cup \Gamma_{out} \cup  \Gamma_{b} \\
	\mbox{and initial conditions.}
\end{array}
\right.
\end{equation}

The possibility to play with the STE term, we allow to simulate different size of problems, different scales in space (small scales with STE, big scales without). To each case, we want determine useful Absorbing Boundary Conditions to truncate space to a more reduce part (the domain $\Omega$). \\

We want to do that by staying close to the transparent operator, the Dirichlet to Neumann map with the non-local properties is not efficient for computations. So we choose to get an approximation of this perfect operator, as proposed by Engquist and Majda \cite{Engquist-Majda-1977}, here with $f_N(X)$ the $(2N+1)$th-order Padé approximation of the square root "$f(X) = \sqrt{1 + X}$" classically written as the rational function \cite{Modave-Warburton}:
\begin{equation}\label{Padetransform2}
	\displaystyle{f_{2N+1}^{\mbox{\scriptsize Padé}}(X) = 1 + \frac{2}{M} \sum_{n=1}^{N} c_n (1 - \frac{1 + c_n}{1 + c_n + X})}
\end{equation}
with ($M=2N+1$):
\begin{equation} \label{Pade-an-bn-cn}
	\displaystyle{c_n = tan^2 (n \pi/M).}
\end{equation}

\section{Classical Padé-Type HABC with the wave equation (a summary)} \label{base}

Before the presentation of the determination of the Padé-type HABC for the CHWM, we summarize the used method on a classical problem to establish the principles, the way, the more important notations. \\

Here we get an approximate non reflective boundary condition with the exact non reflective boundary condition \cite{Engquist-Majda-1977} associated to the classical wave equation:
\begin{equation} \label{WE}
	\partial_t^2 \varphi - c^2 \Delta \varphi = 0 \mbox{ in } \R^2,
\end{equation}
or:
\begin{equation} \label{WEbeforesquareroot}
	c^2 \partial_x^2 \varphi = (\partial_t^2 - c^2 (\Delta - \partial_x^2)) \varphi \mbox{ in } \R^2.
\end{equation}
We consider only the domain $\Omega^- = \{\mathbf{x} \in \R^2 : x < 0\}$ ($x$ the abscissa of the point $\mathbf{x}$) and the perfect boundary condition applied to $\Gamma = \{\mathbf{x} \in \R^2 : x = 0\}$ is (with the square root applied to (\ref{WEbeforesquareroot})) as presented more completely in \cite{Modave-Warburton}:
\begin{equation} \label{WEexactnonreflective}
	c ~ \partial_x \varphi = - \left(\partial_t \sqrt{1 - c^2 (\Delta - \partial_x^2 )/\partial_t^2}\right) \varphi \mbox{ on } \Gamma.
\end{equation}
This perfect non reflective operator associated to (\ref{WEexactnonreflective}) is non local in time and space with the square root in the Dirichlet to Neumann operator. We must approximate the perfect operator to numerical computations. In the case of the propagative waves, the classical way uses the Padé approximations (\ref{Padetransform2}) (more stable than with the Taylor series) to the square root $\sqrt{1+X}$ with:
\begin{equation}
	X=- c^2 (\Delta - \partial_x^2 )/\partial_t^2 .
\end{equation}
So, we get classically:
\begin{equation} \label{WEapproxnonreflective}
	\partial_t \varphi + c ~ \partial_x \varphi = \frac{2}{M} \sum_{n=1}^N c_n \partial_t \left( \varphi_n - \varphi \right) \mbox{ on } \Gamma.
\end{equation}
with each auxiliary fields $\varphi_n$ verifying :
\begin{equation}
	\forall n = 1, N: (1+c_n) \partial_t^2 (\varphi_n - \varphi) - c^2(\Delta - \partial_x^2) \varphi_n = 0 \mbox{ on } \Gamma.
\end{equation}
The auxiliary variables are defined as following (see (\ref{Padetransform2})) with initial conditions equal to zero:
\begin{equation} \label{def_aux_var}
	\varphi_n = \frac{1 + c_n}{1 + c_n + X} ~ \varphi.
\end{equation}
By linearity, the variables $\varphi_n$ are also governed by the same problem as $\varphi$.\\

For our CHWM (\ref{HWMeq}), we apply the same methodology for our specific square root function.

\section{Coupled Padé-type HABC with the CHWM}

In this part, we want to get the Padé-type HABC associated to the complete CHWM (\ref{HWMeq}). We write the coupled models without boundary conditions except for the coupling condition on $\Gamma_s$ (with $f_s$ equal to zero, close to $P_{in}$ and $P_{out}$):
\begin{equation} \label{complete_HWM}
	\left\lbrace
	\begin{array}{ll}
		\displaystyle{ \varepsilon \rho \partial_t^2 \eta - \sigma \partial_x^2 \eta + \rho g \eta + \rho \partial_t \varphi = 0 \mbox{ on } \Gamma_s,} & (a)\\~\\
		\displaystyle{ \partial_t^2 \varphi - c_f^2 \Delta \varphi = 0 \mbox{ in } \Omega, \partial_{\nu} \varphi = \partial_t \eta \mbox{ on } \Gamma_s}. & (b)
	\end{array}
	\right.
\end{equation}
To get Padé-type HABC, we begin to write the transparent operator associated to the basin model. For that, we replicate the previous methodology (see section \ref{base}).
With the way proposed by Stoker \cite{Stoker}, we derive in time the equation (\ref{complete_HWM}.a), we get formally, with the operator $\mathcal{A}$:
\begin{equation}
	\mathcal{A} = \varepsilon \rho \partial_t^2 - \sigma \partial_x^2 + \rho g,
\end{equation}
the time derivative:
\begin{equation}
	\partial_t \eta = - \rho \mathcal{A}^{-1} \partial_t^2 \varphi.
\end{equation}
The basin model (\ref{complete_HWM}.b) can become:
\begin{equation} \label{varphi_model1}
	\left\{
	\begin{array}{l}
		\partial_t^2 \varphi - c_f^2 \Delta \varphi = 0 \mbox{ in } \Omega, \\
		\displaystyle{\partial_{\nu} \varphi = - \rho \mathcal{A}^{-1} \partial_t^2 \varphi \mbox{ on } \Gamma_s.}
	\end{array}
	\right.
\end{equation}
So with the use of the variational formulation, we get ($\forall v \in H_0^1(\Omega)$):
\begin{equation}
	\int_{\Omega} \partial_t^2 \varphi v + \int_{\Omega} c_f^2 \nabla \varphi . \nabla v + \int_{\Gamma_s} c_f^2 \rho \mathcal{A}^{-1} \partial_t^2 \varphi v - \int_{\partial \Omega/\Gamma_s} c_f^2 \partial_{\nu} \varphi v = 0,
\end{equation}
with the dirac operator on $\Gamma_s$ defined by ($\mathbf{x} \in \R^2$):
\begin{equation}
	\int_{\Omega} \delta_{\Gamma_s} f(\mathbf{x}) d\Omega = \int_{\Gamma_s} f(\mathbf{x}) d\Gamma,
\end{equation}
we prefer to write:
\begin{equation}
	\int_{\Omega} ( 1 + c_f^2 \rho A^{-1}  \delta_{\Gamma_s}) \partial_t^2 \varphi v + \int_{\Omega} c_f^2 \nabla \varphi . \nabla v - \int_{\partial \Omega/\Gamma_s} c_f^2 \partial_{\nu} \varphi v = 0,
\end{equation}
A new strong formulation associated can be written:
\begin{equation} \label{BasinStrongFormulationWithDirac}
	\displaystyle{(1 + c_f^2 \rho \mathcal{A}^{-1}  \delta_{\Gamma_s}) \partial_t^2 \varphi - c_f^2 \Delta \varphi = 0 \mbox{ in } \Omega.}
\end{equation}
The last formulation allows to get the Dirichlet to Neumann operator associated to the exact non reflective boundary condition. So we get an equivalent boundary condition as (\ref{WEexactnonreflective}):
\begin{equation} \label{HWM1exactnonreflective}
	c_f ~ \partial_{\nu} \varphi = - \left(\partial_t \sqrt{1 + X}\right) \varphi \mbox{ on } \Gamma
	\mbox{ with } \displaystyle{X = c_f^2 \rho \mathcal{A}^{-1}  \delta_{\Gamma_s} - c_f^2 (\Delta - \partial_x^2)/\partial_t^2}.
\end{equation}
In the case of the propagative waves, we use the Padé-type HABC on $\Gamma$ (here either $\Gamma_{in}$ or $\Gamma_{out}$, Figure \ref{fig:domain}) as (\ref{WEapproxnonreflective}) with  the auxiliary fields $\varphi_n$ ($\forall n = 1,N$) governed by this strong formulation:
\begin{equation}
	(1+c_n) \partial_t^2 (\varphi_n - \varphi) + c_f^2 \rho \mathcal{A}^{-1}  \delta_{\Gamma_s} \partial_t^2 \varphi_n - c^2(\Delta - \partial_x^2) \varphi_n = 0 \mbox{ on } \Gamma.
\end{equation}
For the next steps, we prefer to work with the variational formulation of this last equation ($\forall w \in H_0^1(\Gamma)$):
\begin{equation}
	\int_{\Gamma} (1+c_n) \partial_t^2 (\varphi_n - \varphi) ~ w ~ d\Gamma + \left[ c_f^2 \rho \mathcal{A}^{-1} \partial_t^2 \varphi_n ~ w ~ \right]_{P} + \int_{\Gamma} c^2 \partial_y \varphi_n \partial_y w ~ d\Gamma = 0 \mbox{ on } \Gamma,
\end{equation}
with $P$ equal to $P_{in}$ or $P_{out}$ respectively if $\Gamma$ is equal to $\Gamma_{in}$ or $\Gamma_{out}$ (the other boundary point at the bottom will be resolve by the corner condition as in \cite{Modave-Warburton} or the homogeneous Neumann condition). We take advantage of the last form to get the other strong formulation ($\forall n=1,N$):
\begin{equation} \label{varphi_n_eq1}
	\left\{
	\begin{array}{l}
		(1+c_n) \partial_t^2 (\varphi_n - \varphi) - c_f^2 \partial_y^2 \varphi_n = 0 \mbox{ on } \Gamma,\\
		\displaystyle{\partial_{\nu} \varphi_n = - \rho \mathcal{A}^{-1} \partial_t^2 \varphi_n \mbox{ on } P,}
	\end{array}
	\right.
\end{equation}
with $\varphi$ governed by (\ref{varphi_model1}) or the coupled problem (\ref{complete_HWM}). We can choose to write (\ref{varphi_n_eq1}) also as a coupled problem with the new variables $\eta_n$ ($\forall n=1,N$):
\begin{equation} \label{varphi_eta_n}
	\left\lbrace
	\begin{array}{ll}
		\displaystyle{ \varepsilon \rho \partial_t^2 \eta_n - \sigma \partial_x^2 \eta_n + \rho g \eta_n  + \rho \partial_t \varphi_n = 0 \mbox{ on } \Gamma_s,}\\~\\
		\displaystyle{ (1+c_n) \partial_t^2 (\varphi_n - \varphi) - c_f^2 \partial_y^2 \varphi_n = 0 \mbox{ on } \Gamma, \partial_{\nu} \varphi_n = \partial_t \eta_n \mbox{ on } P}.
	\end{array}
	\right.
\end{equation}
At this step, we separate the study in two parts without and with the STE.

\subsection{Padé-type HABC without surface tension}

For the CHWM without the STE, the operator of the surface model is local. We don't need to use a boundary condition at the $\Gamma_s$ extrema. Directly, we can use the Padé type HABC, first of all with a homogeneous Neumann boundary condition on $\Gamma_b$ to maintain readability. The CHWM without STE can be written with the Padé-type HABC:
\begin{equation} \label{HWMABC0}
	\left\lbrace
	\begin{array}{ll}
		\displaystyle{ \varepsilon \rho \partial_t^2 \eta + \rho g \eta  + \rho \partial_t \varphi = 0 \mbox{ on } \Gamma_s,}\\~\\
		\displaystyle{ \partial_t^2 \varphi - c_f^2 \Delta \varphi = 0 \mbox{ in } \Omega, \partial_{\nu} \varphi = \partial_t \eta \mbox{ on } \Gamma_s,
		\partial_{\nu} \varphi = 0 \mbox{ on } \Gamma_b,} \\~\\
		\mbox{for } \Gamma = \Gamma_{in}, \Gamma_{out} \mbox{ respectively } P = P_{in}, P_{out} \mbox{ and } Q = Q_{in}, Q_{out}: \\
		~
		\begin{array}{ll}
		\displaystyle{ (\partial_t + c_f ~ \partial_{\nu}) \varphi = \frac{2}{M} \sum_{n=1}^N c_n \partial_t \left( \varphi_n - \varphi \right) \mbox{ on } \Gamma,}\\
		\forall n=1,N:\\~\\
			~
			\begin{array}{ll}
			\displaystyle{ \varepsilon \rho \partial_t^2 \eta_n + \rho g \eta_n  + \rho \partial_t \varphi_n = 0 \mbox{ on } P,} \\~\\
			\displaystyle{ \partial_t^2 \varphi_n - c_f^2 \partial_y^2 \varphi_n = 0 \mbox{ on } \Gamma,
			\partial_{\nu} \varphi_n = \partial_t \eta_n \mbox{ on } P,
			\partial_{\nu} \varphi_n = 0 \mbox{ on } Q.}
			\end{array}
		\end{array}
	\end{array}
	\right.
\end{equation}
With or without the inertial term in the surface model, the operator is local.\\

For our numerical computations, we prefer to apply the Padé-type HABC also on $\Gamma_b$. So we use the compatibility conditions for the corners $Q_{in}$ and $Q_{out}$ exactly as \cite{Modave-Warburton}, in $\Omega$ the model is the wave equation, the same model. We modify $(\ref{HWMABC0})$, instead of the homogeneous Neumann boundary conditions on $\Gamma_b$, $Q_{in}$ and $Q_{out}$, we write with the new auxiliary variables $\varphi_m^b (\forall m=1,N)$ on $\Gamma_b$ (with axis given in Figure \ref{fig:domain}):
\begin{equation} \label{HWMABC0-corners}
	\left\lbrace
	\begin{array}{ll}
		\displaystyle{ (\partial_t + c_f ~ \partial_{\nu}) \varphi = \frac{2}{M} \sum_{m=1}^N c_n \partial_t \left( \varphi_m^b - \varphi \right) \mbox{ on } \Gamma_b,}\\~\\
		\displaystyle{ \partial_t^2 \varphi_m^b - c_f^2 \partial_y^2 \varphi_m^b = 0 \mbox{ on } \Gamma_b, \forall m=1,N,}\\~\\
		\mbox{for corners } Q_{in}, Q_{out} , \mbox{ for example with } Q_{out}:\\~\\
		~
		\begin{array}{ll}
		\displaystyle{ (\partial_t + c_f ~ \partial_x) \varphi_m^b = \frac{2}{M} \sum_{n=1}^N c_n ~ \partial_t \frac{(1+c_m) \varphi_n  - c_m \varphi_m^b}{1+c_n+c_m},}\\
		\displaystyle{ (\partial_t - c_f ~ \partial_y) \varphi_n = \frac{2}{M} \sum_{m=1}^N c_m ~ \partial_t \frac{(1+c_n) \varphi_m^b  - c_n \varphi_n}{1+c_n+c_m}.}
		\end{array}
	\end{array}
	\right.
\end{equation}

\subsection{Coupled Padé-type HABC with the surface tension effect}

For the CHWM with the STE, the surface model operator is not local. We need to use a boundary condition at the $\Gamma_s$ extrema.\\

With (\ref{varphi_eta_n}), if we want to get the auxiliary variables $\varphi_n$ ($n=1,N$) useful to the Padé-type HABC (\ref{WEapproxnonreflective}), we must get $\varphi$ and all the variables $\eta_n$. We compute $\varphi$ in $\Omega$ therefore also on $\Gamma$ (i.e. $\Gamma_{in}$ or $\Gamma_{out}$). But we need $\eta_n$ on $P$ (either $P_{in}$ or $P_{out}$), for that with (\ref{varphi_eta_n}), we must compute $\eta_n$ on all $\Gamma_s$ with also $\varphi_n$. It's not very interesting. It's too expensive. We want to do that with the help of a boundary condition for the surface model with the STE.\\

We want to determine a Padé-type HABC applied to the surface model. As previous developments, we must write the transparent operator associated to the surface model. For that, we introduce the operator $\mathcal{D}$ such as:
\begin{equation}
	\varphi = \mathcal{D} \eta,
\end{equation}
represents the basin model (\ref{complete_HWM}.b). The surface model becomes:
\begin{equation} \label{SurfaceStrongFormulationWithDirac}
	\displaystyle{ \varepsilon \rho \partial_t^2 \eta - \sigma \partial_x^2 \eta + \rho g \eta + \rho \mathcal{D} \partial_t \eta = 0 \mbox{ on } \Gamma_s,}
\end{equation}
The last formulation allows to get formally the transparent operator as following ($\varepsilon \neq 0$, for $\Gamma$ equal to $\Gamma_{out}$ with $c_s =\sqrt{\sigma/(\varepsilon \rho)}$):
\begin{equation}
	c_s \partial_x \eta = - \partial_t \sqrt{1 + X} \eta  \mbox{ at } \Gamma \cap \Gamma_s \mbox{ with } X = \frac{g}{\varepsilon}/\partial_t^2 + \frac{1}{\varepsilon} \mathcal{D}/\partial_t.
\end{equation}
The Padé approximation of the square root $\sqrt{1+X}$ is with the auxiliary variables $q_n$ ($\forall n=1,N$):
\begin{equation}
	(\partial_t + c_s \partial_x) \eta = \frac{2}{M} \sum_{n=1}^N c_n \partial_t(q_n - \eta) \mbox{ at } \Gamma \cap \Gamma_s.
\end{equation}
The $N$ auxiliary variables $q_n$ are governed by:
\begin{equation} \label{qn_eq}
	\varepsilon \rho (1 + c_n) \partial_t^2 (q_n - \eta) + \rho g q_n + \rho \mathcal{D} \partial_t q_n = 0 \mbox{ at } \Gamma \cap \Gamma_s.
\end{equation}
At this step, we introduce new auxiliary variables $Q_n$ ($\forall n=1,N$) such as:
\begin{equation}
	Q_n = \mathcal{D} q_n, \forall n=1,N.
\end{equation}
So we can write (\ref{qn_eq}) as following:
\begin{equation}
	\left\lbrace
	\begin{array}{l}
		\displaystyle{\varepsilon \rho (1 + c_n) \partial_t^2 (q_n - \eta) + \rho g q_n + \rho \partial_t Q_n = 0 \mbox{ at } \Gamma \cap \Gamma_s}, \\~\\
		\displaystyle{ \partial_t^2 Q_n - c_f^2 \Delta Q_n = 0 \mbox{ in } \Omega, ~ \partial_{\nu} Q_n = \partial_t q_n \mbox{ on } \Gamma_s}.
	\end{array}
	\right.
\end{equation}
But this problem is not well posed, the variable $q_n$ are not defined on $\Gamma_s$. We add the same surface model (by linearity the variables $q_n$ are also governed by the same problem as $\eta$ (see (\ref{def_aux_var})):
\begin{equation}
	\displaystyle{ \varepsilon \rho \partial_t^2 q_n - \sigma \partial_x^2 q_n + \rho g q_n  + \rho \partial_t Q_n = 0 \mbox{ on } \Gamma_s.}
\end{equation}

For more readability, we summarize the two partial Padé-type HABC:
\begin{equation} \label{FirstVersionOfTwoPadeHABC}
\left\lbrace
\begin{array}{l}
	\displaystyle{ (\partial_t + c_s \partial_{\nu}) \eta = \frac{2}{M} \sum_{n=1}^{N} c_n \partial_t (q_n - \eta) \mbox{ on } \Gamma \cap \Gamma_s,}\\
	\displaystyle{ (\partial_t + c_f \partial_{\nu}) \varphi = \frac{2}{M} \sum_{n=1}^{N} c_n \partial_t (\varphi_n - \varphi) \mbox{ on } \Gamma,}\\
	\mbox{with :}\\
	\left\{
	\begin{array}{l}
		\displaystyle{\varepsilon \rho (1+c_n) \partial_t^2 (q_n - \eta) + \rho g q_n + \rho \partial_t Q_n = 0 \mbox{ on } \Gamma \cap \Gamma_s,}\\
		\displaystyle{ \varepsilon \rho \partial_t^2 q_n - \sigma \partial_x^2 q_n  + \rho g q_n  + \rho \partial_t Q_n = 0 \mbox{ on } \Gamma_s,}\\
		\displaystyle{\partial_t^2 Q_n - c_f^2 \Delta Q_n = 0 \mbox{ in } \Omega, \partial_{\nu} Q_n = \partial_t q_n \mbox{ on } \Gamma_s.}
	\end{array}
	\right. \\~\\
	\left\{
	\begin{array}{ll}
		(1 + c_n ) \partial_t^2 ( \varphi_n - \varphi ) - c_f^2 \partial_y^2 \varphi_n = 0 \mbox{ on } \Gamma, \\
		\displaystyle{ \varepsilon \rho \partial_t^2 \eta_n - \sigma \partial_x^2 \eta_n + \rho g \eta_n + \rho \partial_t \varphi_n = 0 \mbox{ on } \Gamma_s},\\
		\displaystyle{ \partial_t^2 \varphi_n - c_f^2 \Delta \varphi_n = 0 \mbox{ in } \Omega,
			\displaystyle{\partial_{\nu} \varphi_n = \partial_t \eta_n \mbox{ on } \Gamma_s}.}
	\end{array}
	\right.
\end{array}
\right.
\end{equation}
But it's really too expensive to compute such problem. To avoid that, we choose to impose the equalities:
\begin{equation} \label{proposition1}
	q_n = \eta_n \mbox{ and } Q_n = \varphi_n, \forall n=1,N.
\end{equation}
It's less expensive. But it is not so easy. Most of the time, the associated simulations are unstable. We get similar instabilities as written in \cite{Destuynder-Fabre-2016}. To avoid this difficulty, we must use a compatibility condition as following.\\

We propose to introduce a new coefficient in each model, the surface model and the basin model. We do that to get each transparent operator closer to the idea of the Higdon ABC \cite{Higdon-1986}, i.e. with (illustration with the wave equation (\ref{WE})):
\begin{equation}
	a \partial_t + c \partial_{\nu},
\end{equation}
with the $a$ coefficient allowing to modulate the wave velocity. So with the two coefficients $a_s$ and $a_f$, we rewrite the surface model (\ref{SurfaceStrongFormulationWithDirac}) and the basin model (\ref{BasinStrongFormulationWithDirac}):
\begin{equation} \label{StrongFormulationsWithDiracWitha}
	\left\lbrace
	\begin{array}{l}
		\displaystyle{c_s^2 \partial_x^2 \eta = \left( a_s^2 \partial_t^2 + (1-a_s^2) \partial_t^2 + \frac{g}{\varepsilon} + \frac{1}{\varepsilon} \mathcal{D} \partial_t\right) \eta \mbox{ on } \Gamma_s,}\\~\\
	\displaystyle{c_f^2 \partial_x^2 \varphi = \left( a_f^2 \partial_t^2 + (1-a_f^2) \partial_t^2 + c_f^2 \rho A^{-1}  \delta_{\Gamma_s} \partial_t^2 - c_f^2 (\Delta - \partial_x^2) \right) \varphi \mbox{ in } \Omega.}
	\end{array}
	\right.
\end{equation}
Now with the same previous way, we get the new Padé-type HABCs (for example with $\Gamma$ equal to $\Gamma_{out}$):
\begin{equation} \label{PadeHABCwitha}
	\left\lbrace
	\begin{array}{l}
	\displaystyle{(a_s \partial_t + c_s \partial_x) \eta = \frac{2}{M} \sum_{n=1}^N c_n a_s \partial_t(\eta_n^{a_s} - \eta) \mbox{ at } \Gamma \cap \Gamma_s},\\
	\displaystyle{(a_f \partial_t + c_f \partial_x) \varphi = \frac{2}{M} \sum_{n=1}^N c_n a_f \partial_t(\varphi_n^{a_f} - \varphi) \mbox{ on } \Gamma },
	\end{array}
\right.
\end{equation}
with the auxiliary variables $\eta_n^{a_s}$ and $\varphi_n^{a_f}$ defined by ($n=1,N$):
\begin{equation}
	\displaystyle{\eta_n^{a_s} = \frac{1 + c_n}{1 + c_n + X^{a_s}} ~ \eta}, \mbox{ and }
	\displaystyle{\varphi_n^{a_f} = \frac{1 + c_n}{1 + c_n + X^{a_f}} ~ \varphi,}
\end{equation}
with $X^{a_s}$ and $X^{a_f}$ defined by:
\begin{equation}
	\left\lbrace
	\begin{array}{l}
	\displaystyle{X^{a_s} = \frac{(1-a_s^2)}{a_s^2} + \frac{g}{\varepsilon a_s^2}/\partial_t^2 + \frac{1}{\varepsilon a_s^2} \mathcal{D} / \partial_t}, \\~\\
	\displaystyle{X^{a_f} = \frac{(1-a_f^2)}{a_f^2} + \frac{c_f^2 \rho}{a_f^2} A^{-1}  \delta_{\Gamma_s} - \frac{c_f^2}{a_f^2} (\Delta - \partial_x^2)/\partial_t^2}.
	\end{array}
	\right.
\end{equation}
The CHWM is a coupled model. We must therefore look for problems checked by auxiliary variables. Start by rewriting the surface model with:
\begin{equation}
	\mathcal{L}_1(\eta) = \varepsilon \rho \partial_t^2 \eta - \sigma \partial_x^2 \eta + \frac{g}{\varepsilon} \eta \mbox{ and  }
	\mathcal{L}_2(\varphi) = \rho \partial_t \varphi,
\end{equation}
so, we get:
\begin{equation} \label{mathcalEQ}
	\mathcal{L}_1(\eta) + \mathcal{L}_2(\varphi) = 0 \mbox{ on } \Gamma_s .
\end{equation}
Apply the new operators $\mathcal{L}_1$ and $\mathcal{L}_2$ into (\ref{PadeHABCwitha}) for the associated boundaries and with the linearity of these operators, we obtain:
\begin{equation}
	\left\lbrace
	\begin{array}{l}
		\displaystyle{(\partial_t + \frac{c_s}{a_s} \partial_x) \mathcal{L}_1(\eta) = \frac{2}{M} \sum_{n=1}^N c_n \partial_t(\mathcal{L}_1(\eta_n^{a_s}) - \mathcal{L}_1(\eta)) \mbox{ at } \Gamma \cap \Gamma_s},\\
		\displaystyle{(\partial_t + \frac{c_f}{a_f} \partial_x) \mathcal{L}_2(\varphi) = \frac{2}{M} \sum_{n=1}^N c_n \partial_t(\mathcal{L}_2(\varphi_n^{a_f}) - \mathcal{L}_2(\varphi)) \mbox{ on } \Gamma }.
	\end{array}
	\right.
\end{equation}
And adding together the two previous equations, the following condition is useful:
\begin{equation} \label{constraintOna}
	\displaystyle{\frac{c_s}{a_s} = \frac{c_f}{a_f}},
\end{equation}
to get with (\ref{mathcalEQ}) and using the initial conditions equal to zero for the auxiliary variables:
\begin{equation}
	\displaystyle{\sum_{n=1}^N c_n \left(\mathcal{L}_1(\eta_n^{a_s}) + \mathcal{L}_2(\varphi_n^{a_f}) \right) = 0}.
\end{equation}
This equality (\ref{constraintOna}) is similar to the one associated with the Hagstrom-Warburton HABC \cite{Hagstrom-Mar-or-Givoli-2008} and in \cite{Wilk-Jacques-2021}. So it seems an interesting choice that the auxiliary variables $\eta_n^{a_s}$ and $\varphi_n^{a_f}$ verify the same coupled HWM model than $\eta$ and $\varphi$ as implied by the proposal (\ref{proposition1}). But we must link this proposition with the constraint (\ref{constraintOna}).\\

We get our Coupled Padé-type HABC (CP-HABC, abbreviation given page 2). So we write the two Padé-type HABC for the CHWM with the compatibility condition (\ref{constraintOna}):
\begin{equation} 
	q_n = \eta_n  = \eta_n^{a_s} \mbox{ and } Q_n = \varphi_n = \varphi_n^{a_f}, \forall n=1,N,
\end{equation}
as following:
\begin{equation} \label{CP-HABC-eq}
	\left\lbrace
	\begin{array}{l}
		\displaystyle{(a_s \partial_t + c_s \partial_x) \eta = \frac{2}{M} \sum_{n=1}^N c_n a_s \partial_t(\eta_n^{a_s} - \eta) \mbox{ at } \Gamma \cap \Gamma_s},\\
		\displaystyle{(a_f \partial_t + c_f \partial_x) \varphi = \frac{2}{M} \sum_{n=1}^N c_n a_f \partial_t(\varphi_n^{a_f} - \varphi) \mbox{ on } \Gamma },\\
		\mbox{with:}\\
				\begin{array}{l}
					\mbox{on } \Gamma \cap \Gamma_s:\\
					\displaystyle{a_s^2 \varepsilon \rho (1+c_n) \partial_t^2 (\eta_n^{a_s} - \eta) + (1-a_s^2) \varepsilon \rho \partial_t^2 \eta_n^{a_s} + \rho g \eta_n^{a_s} + \rho \partial_t \varphi_n^{a_f} = 0,}\\~\\
					\displaystyle{a_f^2 (1 + c_n ) \partial_t^2 ( \varphi_n^{a_f} - \varphi ) + (1-a_f^2) \partial_t^2 \varphi_n^{a_f} - c_f^2 \partial_y^2 \varphi_n^{a_f} = 0 \mbox{ on } \Gamma},\\ \partial_{\nu} \varphi_n^{a_f} = \partial_t \eta_n^{a_s} \mbox{ on } \Gamma_s.
				\end{array}
	\end{array}
	\right.
\end{equation}
We can note that the auxiliary variables $\eta_n^{a_s}$ and $\varphi_n^{a_f}$ are computed only on $\Gamma$ (here $\Gamma$ is $\Gamma_{out}$, the right border).

\section{Reduction of Padé-type HABC with a big celerity} \label{Specific_Pade_Reduction}

\begin{figure}[h]
	\centering
	\includegraphics[width=11.9cm]{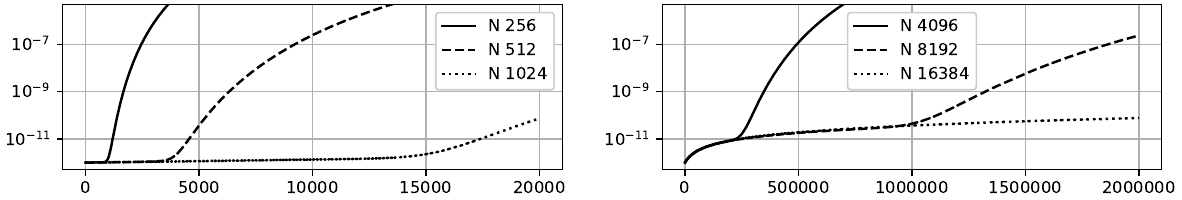}
	\caption{Differences between $\sqrt{1+X}$ and the Padé approximation with several $N$ (given in legends) with $c_f$ equal to $100$ (left graph) and $1000$ $m/s$ (right graph), $X$ (\ref{X_reduc}) in abscissa.}
	\label{fig:bigcf}
\end{figure}

For hydrodynamic problems with a big celerity $c_f$, the operator $X$ (\ref{HWM1exactnonreflective}) becomes (with $\varepsilon$ and $\sigma$ equal to zero to simplify the presentation):
\begin{equation} \label{X_reduc}
	X = c_f^2 (\frac{\delta_{\Gamma_s}}{g} - \partial_{y}^2/\partial_t^2).
\end{equation}
The amplitude of the operator is driven by "$c_f^2$". So with an important value for $c_f$ and in the case of propagative waves, the Padé approximation associated to the square function $\sqrt{1+X}$ must use a large number of terms. Figure \ref{fig:bigcf} illustrates this with different values ($c_f = 100, 1000~m/s$). So the Padé approximation does not seem the better choice.
But with an important celerity coefficient $c_f$ ($c_f \gg 1$), we begin to study the auxiliary variables $\varphi_n$ (first of all with $a_f$ equal to one). The $\varphi_n$ equation on $\Gamma$ with the values $c_n$ (\ref{Pade-an-bn-cn}) can be written:
\begin{equation} \label{Pade-reduc-0}
	\displaystyle{(1 + c_n ) \partial_t^2 (\varphi_n - \varphi) - c_f^2 \partial_y^2 \varphi_n = 0 \mbox{ on } \Gamma},
\end{equation}
At this step, we consider that the time discretizations will be always compatible with the space discretizations, i.e. the ratio between the two discretizations is near or equal to one. In this context, we compare the parameters $(1 + c_n)$ and $c_f^2$. A large part of the $c_n$ values are lower than one and especially an important part of the values $(1 + c_n)$ are lower than $c_f^2$ (for example with $c_f$ equal to $100$ or $1000~m/s$). For sufficiently small $c_n$, the $\varphi_n$ equation can be reduced to:
\begin{equation} \label{Pade-reduc-1}
	\displaystyle{- c_f^2 \partial_y^2 \varphi_n = 0 \mbox{ on } \Gamma}.
\end{equation}
So with the $\varphi_n$ initial conditions equal to zero and the previous simplification, the auxiliary variables $\varphi_n$ remain null. In practical with $c_n$ sufficiently small and an important value for $c_f$, the variables remain near to zero. We propose to neglect these values. \\

\textbf{Remark:} This proposition is enforced with larger time steps for example with implicit schemes, i.e. ratio greater than or equal to one $\geq 0$.

\begin{figure}[h]
	\centering
	\includegraphics[width=4.6cm]{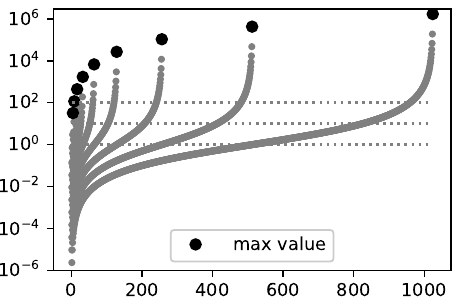}
	\scriptsize{
		\begin{tabular}[b]{|r|rr|rr|rr|}
			\hline
			$n_{abc}$ &  $c_n > 1$ & $\%$ &  $c_n > 10$ & $\%$  &  $c_n > 100$ & $\%$\\	\hline
			4 &  2 &50 &  1 &25  & 0 &0\\	\hline
			8 &  4 &50 &  2 &25  & 1 &12\\	\hline
			16 &  8 &50 &  3 &19  & 1 &6\\	\hline
			32 & 16 &50 &  6 &19  & 2 &6\\	\hline
			64 & 32 &50 & 13 &20  & 4 &6\\	\hline
			128 & 64 &50 & 25 &20  & 8 &6\\	\hline
			256 &128 &50 & 50 &20  &16 &6\\	\hline
			512 &256 &50 &100 &20  &33 &6\\	\hline
			1024 &512 &50 &200 &20  &65 &6\\ \hline
			
	\end{tabular}}
	\caption{Left Graph: $c_n$ values (\ref{Pade-an-bn-cn}) $\forall n=1,N$ for $N = 4, 8, 16, 32, 64, 256, 512, 1024$. The maximum of the set $\{c_n, \forall n=1,N\}$, marked with $\bullet$, grows with $N$. The smaller dot lines gives the values $1, 10, 100$ (the threshold values of the right table) - Right table: Number of the $c_n$ values higher than 1 or 10 or 100 for different orders of the Padé transform.}
	\label{fig:maxcn}
\end{figure}

With Figure \ref{fig:maxcn} (left part), we can note than for small $N$, the maximum value of $c_n$ is small, too small to get a correct approximation with the Padé transform. And of course, it seems useful to use a sufficiently big $N$ to get a bigger $c_N$ and to get more none neglected terms of the auxiliary variables to compute a sufficient Padé approximation of the square root (\ref{WEapproxnonreflective}). Figure \ref{fig:maxcn} (right part) shows also different numbers of $c_n$ coefficients higher than different thresholds.
So with big $c_f$ as $100$ or $1000~m/s$, the thresholds $10$ and $100$ must be interesting to use to reduce the number of the terms of the Padé approximation. As we see in the next section for the big $c_f$ as $1000~m/s$, the threshold can be more important.\\

We must study also the auxiliary variables $\eta_n$  of the Padé-type HABC associated to the surface model. The reduction possibilities are limited by the value of the coefficient $\varepsilon$. We have (with the simplification $a_s$ equal to one):
\begin{eqnarray} \label{eq-aux-eta-v0}
	\displaystyle{\varepsilon \rho (1+c_n) \partial_t^2 (\eta_n - \eta) + \rho g \eta_n + \rho \partial_t \varphi_n = 0 \mbox{ on } \Gamma \cap \Gamma_s}.
\end{eqnarray}
The reduction method is efficient with $\varepsilon$ sufficiently small or without the added mass term, the relation connects directly $\eta_n$ to $\varphi_n$. For sufficiently small $c_n$, the auxiliary variables $\eta_n$ remain near to zero. These terms can be also neglected.\\

With the compatibility condition (\ref{constraintOna}), we need to add some more explanations. We must make a choice with the compatibility condition (\ref{constraintOna}). We discuss two opposite choices:
\begin{eqnarray} \label{reduction-cases}
\left\lbrace
\begin{array}{cc}
a_f = 1, a_s = c_s/c_f,& (a)\\
a_s = 1, a_f = c_f/c_s. & (b)
\end{array}
\right.
\end{eqnarray}
The case (\ref{reduction-cases}.a) changes only the $\eta_n$ equation (\ref{eq-aux-eta-v0}) as in (\ref{CP-HABC-eq}):
\begin{eqnarray}
\displaystyle{a_s^2 \varepsilon \rho (1+c_n) \partial_t^2 (\eta_n^{a_s} - \eta) + (1-a_s^2) \varepsilon \rho \partial_t^2 \eta_n^{a_s} + \rho g \eta_n^{a_s} + \rho \partial_t \varphi_n^{a_f} = 0}.
\end{eqnarray}
For the water case and $\varepsilon \in [10^{-9}, 10^{-3}]$, we get $a_s$ lower than $1$. So we don't loose the previous reduction possibilities, the negligible terms are always negligible. \\
But in the case (\ref{reduction-cases}.b), the $\varphi_n$ equation (\ref{Pade-reduc-0}) changes as in (\ref{CP-HABC-eq}):
\begin{eqnarray} \label{eq-aux-varphi-v1}
\displaystyle{a_f^2 (1 + c_n ) \partial_t^2 ( \varphi_n^{a_f} - \varphi ) + (1-a_f^2) \partial_t^2 \varphi_n^{a_f} - c_f^2 \partial_y^2 \varphi_n^{a_f} = 0}.
\end{eqnarray}
For the water case and $\varepsilon \in [10^{-9}, 10^{-3}]$, we get $a_f$ upper than $1$. So for example with $a_f^2 \gg 1$, the value "$a_f^2$" can be closer than $c_f^2$ and also $a_f^2 c_n$. The two first terms of (\ref{eq-aux-varphi-v1}) are no longer as negligible as before. The reduction is less efficient. So we prefer to use the case (\ref{reduction-cases}.a).\\

For problems with an important celerity coefficient $c_f$ ($c_f \gg 1$) and a reasonable added mass on the surface, it's useful to use a sufficiently big $N$ and we can/must avoid to compute all the auxiliary variables for the sufficiently small $c_n$ with the case (\ref{reduction-cases}.a) and compatible space $h$ and time $\Delta t$ discretizations (in particular with $\Delta t \geq h$). We use this statement to reduce efficiently the size, the computational cost of our numerical problems.

\section{Numerical simulations}

We solve the problem (\ref{HWMeq}) associated to the CHWM without and with the STE. On the boundaries $\Gamma_{in}$, $\Gamma_{out}$ and $\Gamma_b$, we apply ours HABCs. Without STE, we use (\ref{HWMABC0}) with the compatibility conditions for the corners $Q_{in}$ and $Q_{out}$ (\ref{HWMABC0-corners}) exactly as \cite{Modave-Warburton}. With the STE, we use our CP-HABC (\ref{CP-HABC-eq}) with the constraint (\ref{constraintOna}) and the compatibility conditions for the corners (\ref{HWMABC0-corners}).

\subsection{Space and time discretizations}
For the space discretization of the domain, we use the Finite Element formulation with an uniform mesh with parallelepiped elements. We use the standard spatial Galerkin Finite Element discretization with the polynomial functions of Lagrange of order four. \\
For the time discretization, we choose the implicit version of the Newmark scheme ($\gamma = 0.5, \beta = 0.25$) not here, in this article, to use different time steps, for example designed for fast and slow domains, but we use in almost all our applications the same time step designed for surface waves, the scheme in the implicit version allowing an appreciate robustness. It's rather the change of space discretization that affects the quality of results associated with the use of HABC. Just for one case (a special case with an object below the surface), we compute with different time steps\\
The space and time discretizations are chosen to describe the behavior associated to the surface waves. The used space, time steps respectively $h$ and $\Delta t$ are given in the caption of Table \ref{table-cas} with the physical characteristics of the different cases.

\subsection{Different cases}
\begin{figure}[h!]
	\centering
	\includegraphics[width=10cm]{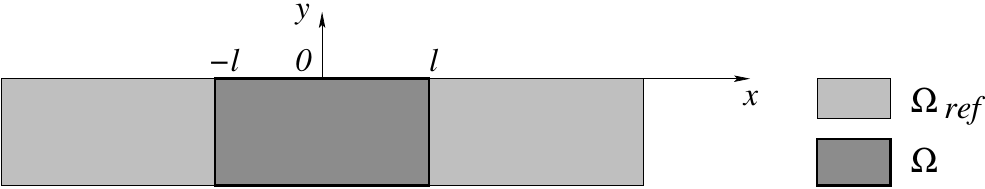}
	\caption{The computational domains (see also Figure \ref{fig:domain}).}
	\label{fig:domain-appli}
\end{figure}

The truncated domain $\Omega$ is $[-l, l]\times[-l,0]$ (Figure \ref{fig:domain-appli}). 
The reference domain $\Omega_{ref}$ is $[-l_{ref}, l_{ref}]\times[-l,0]$. The two domains $\Omega$ and $\Omega_{ref}$ are at the same depth.
The parameter $l_{ref}$ and the final time $T$ depend on the values of $\rho$, $g$, $\sigma$ and $l$. The classical hydrodynamic model uses the water parameters:
\begin{equation} \label{water parameters}
	\rho = 1000 ~Kg.m^{-3}, \sigma = 0.075 ~N.m^{-1}, c_f = 1500 ~m.s^{-1}.
\end{equation}
We want to test this fluid in the gravity $g$ equal to $9.81 ~m.s^{-2}$ and especially for medium and small scales in space. We want to test also this model for other materials, for example as Helium liquid (4°K). In this last case, $c_f$ can be equal to $200 ~m.s^{-1}$ (sound 1) or $20 ~m.s^{-1}$ (sound 2, heat waves, helium II with $\sigma$ hundred to thousand times smaller \cite{Guyon}, \cite{Caupin}).
But in this publication to avoid too many cases (computations), we limit our tests with the parameters $\rho$ and $\sigma$ as (\ref{water parameters}) and:
\begin{equation}
	c_f = 1, 100, 1000 ~m.s^{-1} \mbox{ with } g = 9.81~m.s^{-2},
\end{equation}
and just one case with $\sigma$ equal to zero. In particular, we use $c_f$ equal to one to get some various references (also with $\sigma$ equal to zero) useful to compare the precision with the other more real cases with the specific reduction of the Padé coefficients (Section \ref{Specific_Pade_Reduction}). Table \ref{table-cas} gives the values of the parameters of the different cases of our numerical simulations. For each case, we use three meshes, the coarser mesh with an element size equal to $h$, a better with $h/2$ and the best with $h/4$ ($h$ given in the caption of Table \ref{table-cas}).\\

\begin{table}[h!]
	\begin{tabular}{|c|c|c|c|c|c|c|c|c|c|c|} \hline
			case   & $c_f$ & $\sigma$ & $T$ & $\varepsilon$ & $~~~l, l_{ref}$\\
		     & \scriptsize{$m/s$} & \scriptsize{$N/m$} &  \scriptsize{$s$} &    &  \scriptsize{$m, m$} \\ \hline
			1 & 1.  &  0. &  1.5  &  0  & 0.1, 1.0\\ \hline 
			\hline
			11 & 1.  &  0.075 &  0.9  &  1.e-3  & 0.1, 0.5\\ \hline
			12 & 1.  &  0.075 &  0.9  &  1.e-9  & 0.1, 0.5\\ \hline
			\hline
			\multicolumn{6}{|c|}{big $c_f$, $c_n$ reduction} \\ \hline
			211 & 100.  &  0.075 &  0.9  &  1.e-3  & 0.1, 0.5\\ \hline
			311 & 1000.  &  0.075 &  0.9  &  1.e-3  & 0.1, none\\ \hline
	\end{tabular}
	\caption{The different computational cases with $h=0.0025~m$ with $\Delta t = 0.002~s$.} \label{table-cas}
\end{table}
With the initial conditions equal to zero, the excitation is given by $f_s$ (\ref{HWMeq}) (with $T_{e} = T_{excit}= 0.1 s, A = 1000, n_f = 20, x_0 = 0$) useful to test with different time and space frequencies:
\begin{equation} \label{excitation}
	\begin{array}{r|cl}
		f_s(x,t) =& \displaystyle{\sum_{n=1}^{n_f} A e^{-10 (\frac{n-1}{n_f})^2} sin(\frac{2n\pi t}{T_{e}}) ~ e^{-100(\frac{x-x_0}{l})^2}} & \mbox{if } t \le T_{excit},\\
		& 0 & \mbox{if } t > T_{excit}.
	\end{array}
\end{equation}

\begin{figure}[h]
	\centering
	\includegraphics[width=5.5cm]{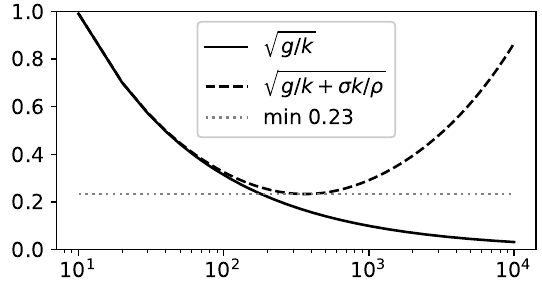}
	\caption{Phase speed $C_p(k)$ of the surface waves $\omega/k$ with and without STE for $g= 9.81 m/s^2$ and $\varepsilon = 0$. The values ($\rho = 1000 Kg/m^3, g = 9.81 m/s^2, \sigma = 0.075 N/m$) are near of the water parameters.}
	\label{fig:phasespeed}
\end{figure}

We choose the final time $T$ (Table \ref{table-cas}) to allow a large part of the initial waves to go out of the domain $\Omega$ and particulary for the waves on the surface. For that we can use the phase speed associated to the surface waves (Figure \ref{fig:phasespeed}). \\

For the reference problems with $c_f$ equal to $1$ $m/s$, the value of $l_{ref}$ is sufficient to avoid the return of the waves in the domain $\Omega$. For $c_f$ equal to $100$ and $1000$ $m/s$, the correct lengths to avoid the return of the faster waves are not possible for our computer resources. So we use on the output boundaries our CP-HABCs with a sufficient high order. Nevertheless, the length $l_{ref}$ is sized with the phase speed of the waves on the surface to avoid the return of the waves associated to the surface waves.\\

For the compatibility condition(\ref{constraintOna}), we use the efficient case (\ref{reduction-cases}.a), i.e. $a_f$ equal to one. This condition is used especially with a big celerity as for water. This compatibility condition is not so useful for the cases with $c_f$ equal to $1~m.s^{-1}$.\\

To analyze the precision of the numerical simulations, we compute these different errors ($t \in [0,T]$):
\begin{equation} \label{errors}
	\left\lbrace
	\begin{array}{l}
		\displaystyle{e_{\eta}(t) = \frac{||\eta - \eta_{ref}||_{\infty, \Gamma}}{||\eta_{ref}||_{\infty, \Gamma \times [0,T]}} \mbox{ and } E_{\eta} = ||e_{\eta}||_{2, [0,T]}},\\~\\
		\displaystyle{e_{\varphi}(t) = \frac{||\varphi - \varphi_{ref}||_{\infty, \Omega}}{||\varphi_{ref}||_{\infty, \Omega \times [0,T]}} \mbox{ and } E_{\varphi} = ||e_{\varphi}||_{2, [0,T]}},\\
	\end{array}
	\right.
\end{equation}
For the case $c_f$ equal to $100$ $m/s$, we compute also the energy in $\Omega$ and on the surface. And for the case $c_f$ equal to $1000$ $m/s$, we compute only the energy. The previous errors (\ref{errors}) are too expensive to get for the last case.

\subsection{Special test with an object below the surface} \label{special_case}

We propose a last test with a small elliptical shaped object just below the surface with a more important final time ($T=5\times0.9~s$) allowing a longer excitation time and the wave celerity $c_f$ equal to $1500~m/s$ (water celerity). The other physical properties are those of the previous case 311. The presence of an object just below the surface is a motivation presented in our introduction, for example interesting for fluid-structure interaction with Lagrangian or/and Eulerian methods with fitted or unfitted mesh approaches and level set method (IBM, FD, ALE \cite{IBM-Peskin}, \cite{FD-Glowinski}, \cite{ALE-Donea}, \cite{levelset-Cottet}). The longer excitation time is interesting to evaluate the robustness of the method in this context. \\

The domain $\Omega$ is $[-l, l] \times [-l/4, 0]$ ($l=0.1~m$ as previous cases) with the same space discretization and the same time discretization. The center of the elliptical shaped object is $(0.05, -0.01)$, the two characteristic lengths are $0.01~m$ and $0.005~m$. The excitation $f_s$ is always given by (\ref{excitation}) but with a double $T_e$ equal to $0.2~s$, $T_{excit}$ equal to $12~T_e$ and $x_0$ equal to $-l/2$. With this special case, we use three different times steps ($\Delta t, \Delta t/2, \Delta t/4$) especially to show the effect on the reduction of the Padé approximation (Section \ref{Specific_Pade_Reduction}).

\subsection{Numerical results - Cases with $c_f$ equal to one}
\begin{figure}[h]
	\centering
	\includegraphics[width=11.9cm]{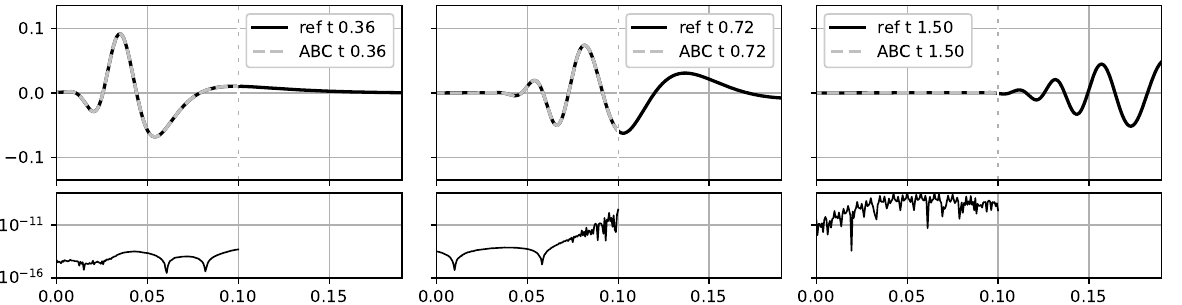}
	\caption{Case 1, ABC order 32 - The top graphs show $\eta$ computed with the Pade-type ABC (dashed line) and the reference (solid line) for three different times given with the legend, the y-axis is the same for all the graphs. The small graphs just below the top graphs show the associated errors $|\eta - \eta_{ref}|$, the y-axis is the same for all the graphs.}
	\label{fig:resumeetacas1clpadeeva0q4nbase20abc32visu1d150}
\end{figure}
\begin{figure}[h]
	\centering
	\includegraphics[width=11.9cm]{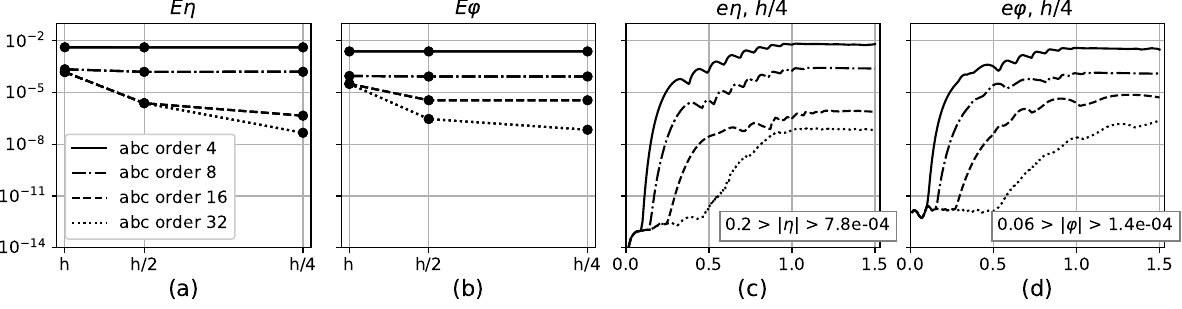}
	\caption{Case 1 - The graphs (a) and (b) show respectively the errors $E_{\eta}$ and $E_{\varphi}$ (\ref{errors}) for different mesh discretizations given in abscissa and for different orders of the CP-HABC given by the legend of the graph (a). The legend and the log10 y-axis of the graph (a) are the same for the other graphs. The graphs (c) and (d) show respectively the errors $e_{\eta}(t)$ and $e_{\varphi}(t)$ (\ref{errors}) (time in abscissa) for different orders of the CP-HABC for the best mesh discretization. In the lower right corner, we can see the maximum and the last ($t=T$) amplitudes for $\eta$ (graph (c)) and $\varphi$ (graph (d)).}
	\label{fig:resumetraceetaphicas142032}
\end{figure}
The first cases $\{1, 11,12\}$ use $c_f$ equal to one with the other parameters near of the water parameters (see Table \ref{table-cas}).\\
 
For the case 1, the only case with $\sigma$ equal to zero, Figure \ref{fig:resumeetacas1clpadeeva0q4nbase20abc32visu1d150} shows the evolution of the amplitudes $\eta$ and $\eta_{ref}$ for three different times. With this presentation, the two solutions $\eta$ and $\eta_{ref}$ merge together, the limit between the two computational domains is marked at $0.1~m$ by a vertical dashed line. At the final time, a large part of the waves are out without a visible difference. The bottom small graphs give the absolute difference $|\eta -\eta_{ref}|$ for the same times. The differences are always very small below $10^{-8}$. Figure \ref{fig:resumetraceetaphicas142032} shows the errors $E_{\eta}$ and $E_{\varphi}$ for different mesh discretizations (\ref{errors}) and for different orders of the CP-HABC, the higher the order of the condition the better the accuracy. The high orders of the CP-HABC need sufficient space discretization levels to work efficiently. The best result is obtained with the finer mesh ($h/4$). Figure \ref{fig:resumetraceetaphicas142032} shows also the errors $e_{\eta}(t)$ and $e_{\varphi}(t)$ (\ref{errors}) for the best/fine mesh and for different orders of the CP-HABC. The CP-HABC performs as expected, the errors decrease with increasing the ABC order to get very worst small values near $10^{-7}$ for the relative errors.\\

\begin{figure}[h]
	\centering
	\includegraphics[width=11.9cm]{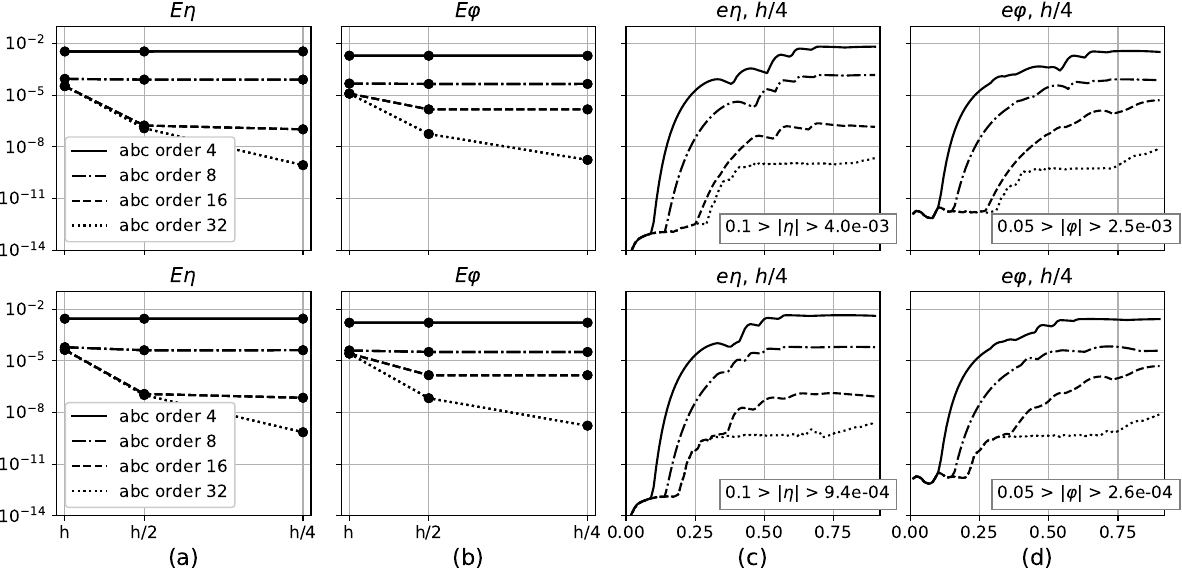}
	\caption{Cases 11, 12 respectively top line, bottom line - Each line uses similar presentation as Figure \ref{fig:resumetraceetaphicas142032}.}
	\label{fig:resumetraceetaphicas1142032}
\end{figure}
\begin{figure}[h]
	\centering	
	\includegraphics[width=11.9cm]{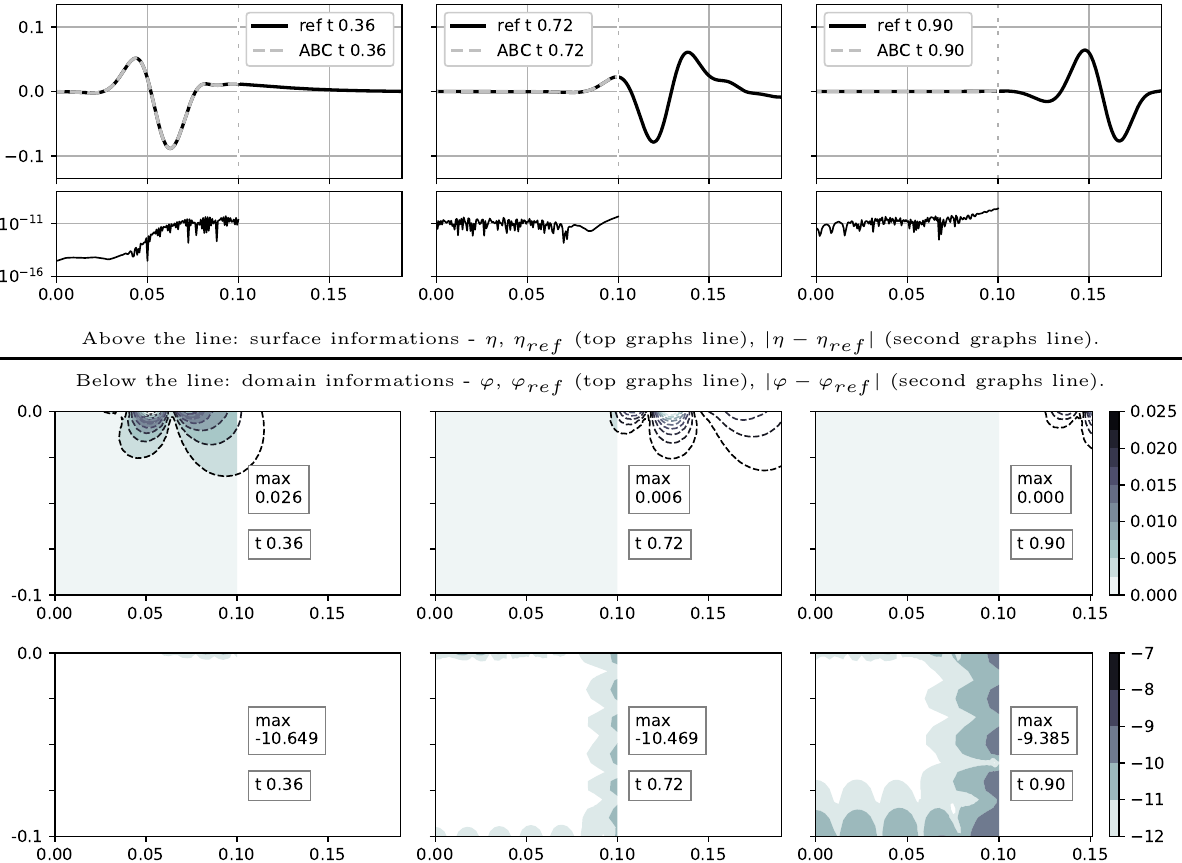}
	\caption{Case 12, ABC order 32 - Above the line, similar presentation as Figure \ref{fig:resumeetacas1clpadeeva0q4nbase20abc32visu1d150}. The next graphs just below the line show $\varphi$  computed with the Pade-type ABC (fill contours with the colorbar on the right side) and the reference (dashed contours) for the same times. The bottom graphs show the log10 of the error $|\varphi - \varphi_{ref}|$, the colorbar gives the associated power. In each bottom graphs, the legend gives the maximum of the log10 of the error $|\varphi - \varphi_{ref}|$ and the time. The left y-axis is the same for all graphs of each lines.}
	\label{fig:resumeetacas11clpadeeva0q4nbase20abc32visu1d90}
\end{figure}
The cases 11 and 12 with the STE, respectively with added mass on the surface and almost not, give similar results. Figure \ref{fig:resumetraceetaphicas1142032} shows the errors. The CP-HABC performs as expected. The behaviors are similar to the previous case. The relative errors are a little smaller than previous case, the smaller values are around $10^{-9}$ for the best cases. For the case 12 with almost no added mass on the surface, Figure \ref{fig:resumeetacas11clpadeeva0q4nbase20abc32visu1d90} allows to show the output of the waves for $\eta$ and $\varphi$ with $|\eta - \eta_{ref}|$ on $\Gamma_s$ and $|\varphi - \varphi_{ref}|$ in $\Omega$ for the best ABC order and the best mesh. The absolute differences $|\eta - \eta_{ref}|$ on $\Gamma_s$ and $|\varphi - \varphi_{ref}|$ in $\Omega$ confirms the interesting quality of the CP-HABC. The worst absolute errors are around $10^{-9}$ with the maximum values for $\eta$ and $\varphi$.\\

We can note that the STE changes the behavior of the solution between the case 1 (without STE, Figure \ref{fig:resumeetacas1clpadeeva0q4nbase20abc32visu1d150}) and the case 12 (with STE, Figure \ref{fig:resumeetacas11clpadeeva0q4nbase20abc32visu1d90}). \\

So with theses three previous cases with $c_f$ equal to $1~m/s$, we have already efficient reference numerical results, also useful to compare with more realistic cases with bigger values for $c_f$, as we see in the next section.

\subsection{Numerical results - Cases with big $c_f$ with Padé-type HABC reduction}

\begin{figure}[h]
	\centering
	\includegraphics[width=11.9cm]{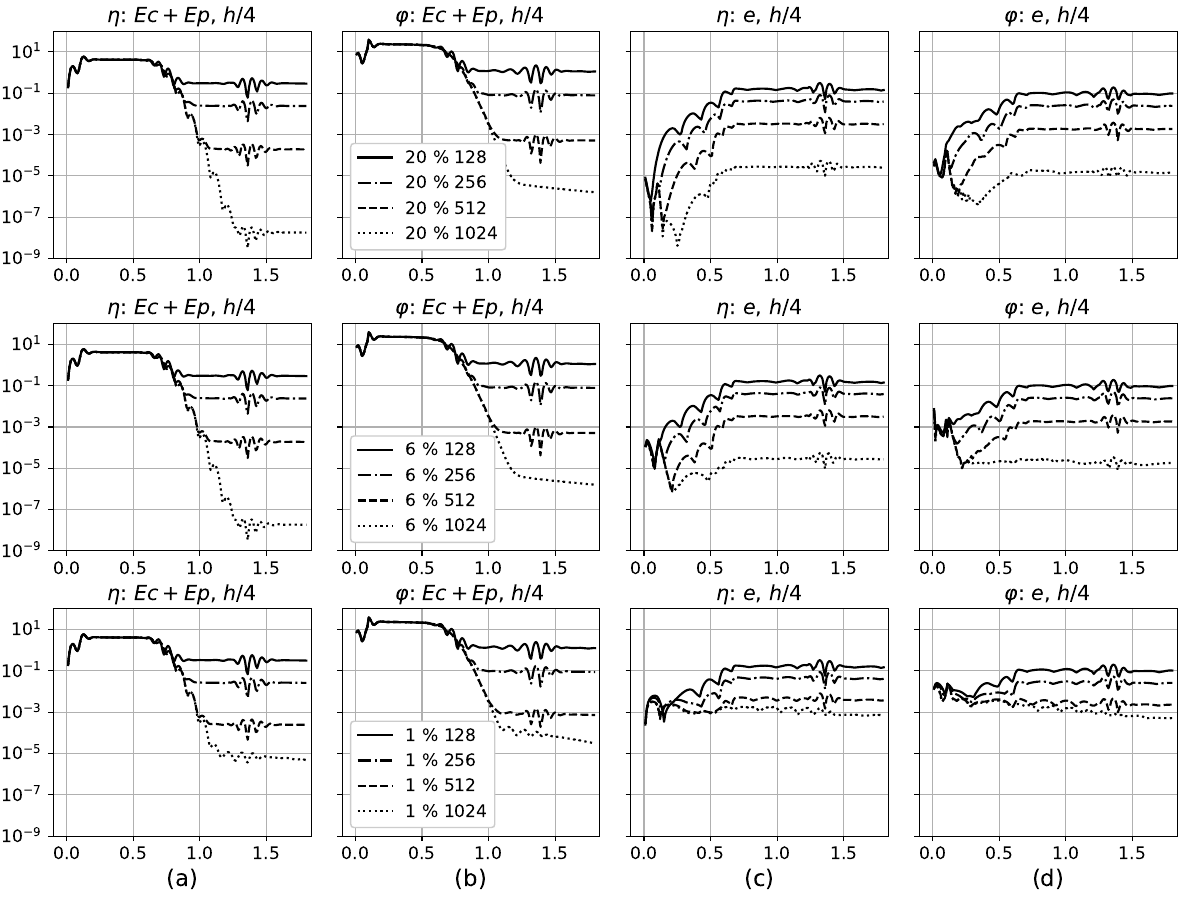}
	\caption{Case 211 with different reductions levels associated to the Padé transform, we keep $20 \%$ of the last auxiliary variables (top line), $6 \%$ (middle line), $1 \%$ (bottom line) and for different orders (128, 256, 512, 1024 marked with the (b) column legend) - Function of time, the (a) column gives the total energy for $\eta$ in $\Gamma_s$, the (b) column the total energy for $\varphi$ in $\Omega$, the (c) column the relative errors $e_{\eta}(t)$, the (d) column $e_{\varphi}(t)$. The left y-axis is the same for all the graphs.}
	\label{fig:trace3020211}
\end{figure}
The first case is the 211 (Table \ref{table-cas}) equivalent to the case $11$ but with $c_f$ equal to $100$ $m/s$ and only with the finer mesh ($h/4$). \\

The graph columns (a) and (b) of Figure \ref{fig:trace3020211} show the time evolution of the total energy respectively for $\eta$ on the surface $\Gamma_s$ and for $\varphi$ in the domain $\Omega$. This type of visualization is interesting because it allows to avoid to compute a reference case to evaluate the interest of our CP-HABC associated to our CHWM as we do for the previous cases. The energy allows to get a synthesis of the presence of the waves in the domain $\Omega$ and the surface $\Gamma_s$. Here, we point out again that it's very expensive to use a sufficient length $l_{ref}$ to prevent the return of the faster waves in the domain $\Omega$ (without ABC or HABC on $\Gamma_{in}$ and $\Gamma_{out}$ boundaries of $\Omega_{ref}$). Nevertheless for this case, we compute a reference case but with a sized length $l_{ref}$ sized to avoid the return of the surface waves and with our CP-HABC on the output boundaries $\Gamma_{in}$, $\Gamma_{out}$ and $\Gamma_b$ with all the terms of the Padé approximation with the order equal to 1024 (without reduction). So the columns (c) and (d) (Figure \ref{fig:trace3020211}) show the relative errors $e_{\eta}(t)$ and $e_{\varphi}(t)$ (\ref{errors}) for the case with the finer mesh ($h/4$). We compute the total energy and the errors for different orders of the CP-HABC (128, 256, 512 and 1024) and different reduction levels (section \ref{Specific_Pade_Reduction}) $20$, $6$ and $1$ $\%$ as shown in the legends of the graphs of the (b) column. These cases with the two evaluation modes (based on errors and energy) make the junction between the previous cases with the errors (\ref{errors}) and the next cases (with $c_f$ equal to $1000$ $m/s$) with only the time evolution of the total energy. We can remark that the final time $T$ must be a little bit more important to get the end of the reduction of the total energy in particular for the best case, here with the Padé order equal to $1024$.\\
The CP-HABC performs as expected if we use sufficient high order for the Padé transform	 (see also Figure \ref{fig:bigcf}). The errors decrease with increasing the ABC order for the different reduction levels. For this case, a good balance between the HABC efficient and the cost of the computation (depending of the choose of the reduction level) seems to use between $6$ and $1$ $\%$ of the last terms of the Padé approximation. So for the CP-HABC order equal to $1024$, we use only less than between $64$ and $12$ terms. These are still not too big numbers to get efficient accuracy.\\

\begin{figure}[h]
	\centering
	\includegraphics[width=11.9cm]{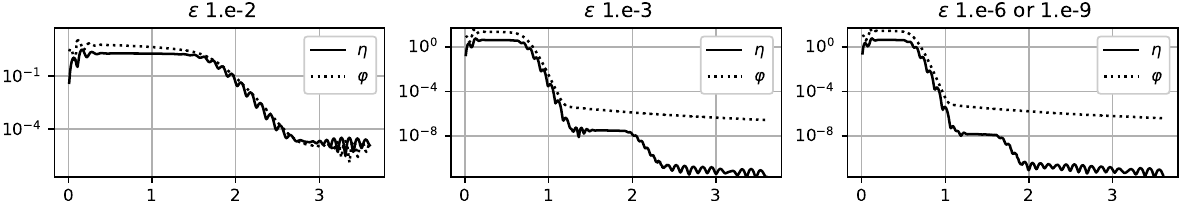}
	\caption{Function of time, the total energy for $\eta$ (solid line) and $\varphi$ (doted line) with different values for $\varepsilon$: $10^{-2}$ (left graph), $10^{-3}$ (middle), $10^{-6}$ or $10^{-9}$ (right) - Same case as Figure \ref{fig:trace3020211} with a reduction of $3 \%$ of the Padé transform for the order equal to 1024 ($T = 3.6~s$).}
	\label{fig:cas3013211r930}
\end{figure}
\textbf{Remark:} As the method allowing to determine the ABC uses the approximation of the Dirichlet to Neumann map. The Padé-Type HABC is always in connection with the continuous model. If the approximation is improved, to get an improvement in accuracy of the numerical simulations, we must use a sufficient fine discretization. For example, an important accuracy is useful to use our CP-HABC to develop an efficient interface condition for Domain Decomposition Method. In direct simulations with just one domain, a very strong precision is not always used. So in this context, we can use a small number of terms (last terms with big N for example with water cases) of the Padé approximation. For example with only $1\%$ of $1024$, we get relative errors around $10^{-3}$.\\

In Figure \ref{fig:trace3020211}, we can note an artefact just before the time equal to $1.5~s$. It's just the very small wave bounces on the opposite boundary coming with a delay associated to the distance between the two boundaries $\Gamma_{in}$ and $\Gamma_{out}$. Nevertheless, these small wave bounces show the loss of precision with the reduction of the Padé-type HABC associated to the surface model. The coefficient $\varepsilon$ here equal to $10^{-3}$ is near the limit, with a bigger $\varepsilon$, for example $10^{-2}$, the errors, the wave bounces are more important as shown in the left graph of Figure \ref{fig:cas3013211r930} around the time $t$ equal to $3~s$ (the log10 y-axis are different between the left graph and the other graphs). With the two others graphs, the time evolutions of the total energies show with a small ($10^{-3}$) and a very small ($10^{-6}$ or $10^{-9}$, the two solutions merge, the added mass becomes negligible) coefficient $\varepsilon$ a good quality. For Figure \ref{fig:cas3013211r930}, we use only $3 \%$ of the last terms of the Padé approximation (i.e. $32$ terms). We can note also the very good behavior for a longer time (here $3.5$ s).\\

\begin{figure}[h]
	\centering
	\includegraphics[width=11.9cm]{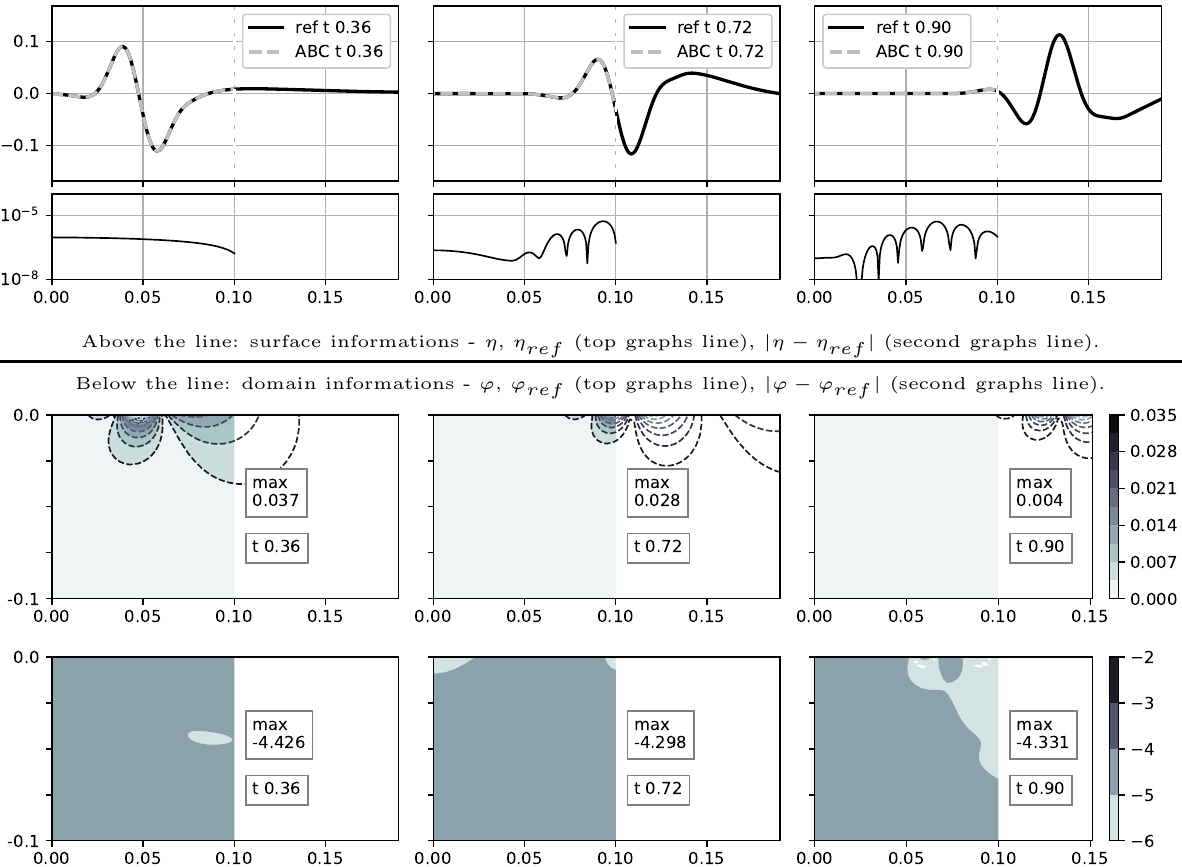}
	\caption{Case 211, ABC order 1024, reduction level $6 \%$ - Similar presentation as Figure \ref{fig:resumeetacas11clpadeeva0q4nbase20abc32visu1d90}.}
	\label{fig:resumeetacas3006211clpadeeva0q4nbase20abc63visu1d90}
\end{figure}
For the case 211 with $6 \%$ reduction level of the Padé approximation terms associated to the order 1024, Figure \ref{fig:resumeetacas3006211clpadeeva0q4nbase20abc63visu1d90} shows for three moments the vertical displacement $\eta$ and the velocity potential $\varphi$ and the absolute errors with the reference solution. We can visualize the correct output of the different waves.\\

\begin{figure}[h]
	\centering
	\includegraphics[width=11.9cm]{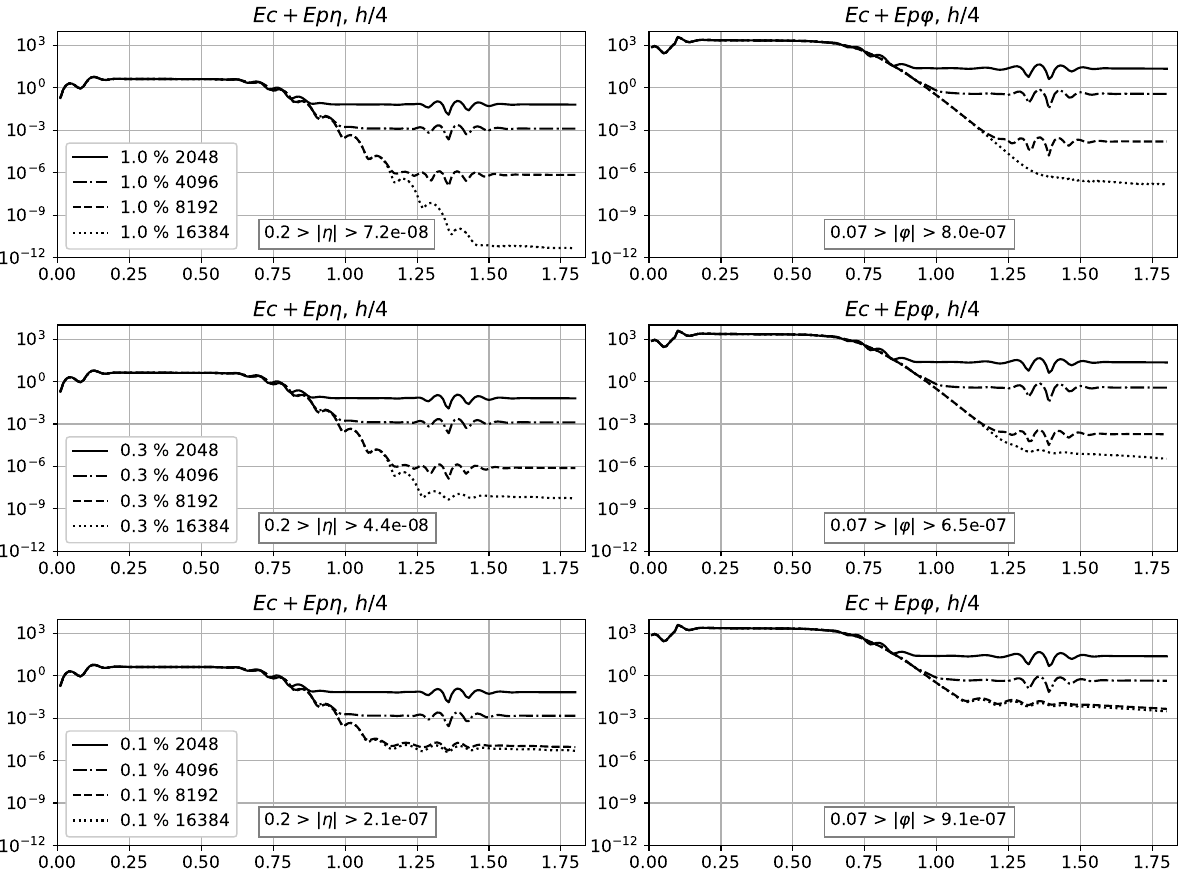}
	\caption{Case 311 ($c_f=1000~m/s$) with different reductions levels associated to the Padé transform, $1 \%$ of the last auxiliary variables (top line), $0.3 \%$ (middle line), $0.1 \%$ (bottom line) and for different orders (2048, 4096, 8192, 16384 marked with the left column legend) - Function of time, the left column gives the total energy for $\eta$ in $\Gamma_s$, the right column the total energy for $\varphi$ in $\Omega$.}
	\label{fig:trace6001211}
\end{figure}
Figure \ref{fig:trace6001211} shows the evolution of the total energy of the case 311 with $c_f$ equal to $1000~m/s$ with the finer mesh ($h/4$) for $\eta$ on the surface $\Gamma_s$ (left column) and for $\varphi$ in the domain $\Omega$ (right column). We compute the total energy for different orders of the CP-HABC with different reduction levels (see section \ref{Specific_Pade_Reduction}). In this case and also the previous case, the coefficient $c_f$ is big ($1000$ and $100~m/s$), we must use very high order to get sufficiently big values for the $c_n$ coefficients to transfert sufficiently $\varphi$ information to $\varphi_n$ with (\ref{Pade-reduc-0}). If the $cn$ coefficient is too small, the $\varphi_n$ equation can be reduced to (\ref{Pade-reduc-1}), so $\varphi_n$ remains null (with null initial conditions and with $\varepsilon$ not so big, see last part of Section \ref{Specific_Pade_Reduction})). With the coefficient $c_f$ equal to $100~m/s$, we have seen than an interest order is around 1024 (Figure \ref{fig:bigcf} and Figure \ref{fig:trace3020211}). With the coefficient $c_f$ equal to $1000~m/s$, we must use higher order (see Figure \ref{fig:bigcf}). The order 1024 is too small to get efficient results, it's why we use the order (2048, 4096, 8192, 16384) and we can use just a little last part of the terms of the Padé approximations, the reduction levels are $1$, $0.3$ and $0.1$ $\%$ as shown in the legends of the graphs of the left column. \\
We get similar behaviors as previous cases. The CP-HABC performs as expected if we use sufficient high order for the Padé transform. For example, the case with the order equal to $16384$ with $0.3 \%$ of the last terms of the Padé approximation allows to get a result with a very interesting accuracy (around $10^{-8}$ for the total energy associated to the surface problem and around $10^{-5}$ ($t \simeq 1.25~s$) for the basin problem). For the case, the Padé approximation uses around $50$ terms. It still remains usable.

\subsection{Numerical results - Special test}

For the test with the elliptical object below the surface with a long simulation time and physical properties with STE equivalent to water (see subsection \ref{special_case}), we test our CP-HABC with strong reductions of the Padé approximation and different time steps. \\

\begin{figure}[h]
	\centering
	\includegraphics[width=11.9cm]{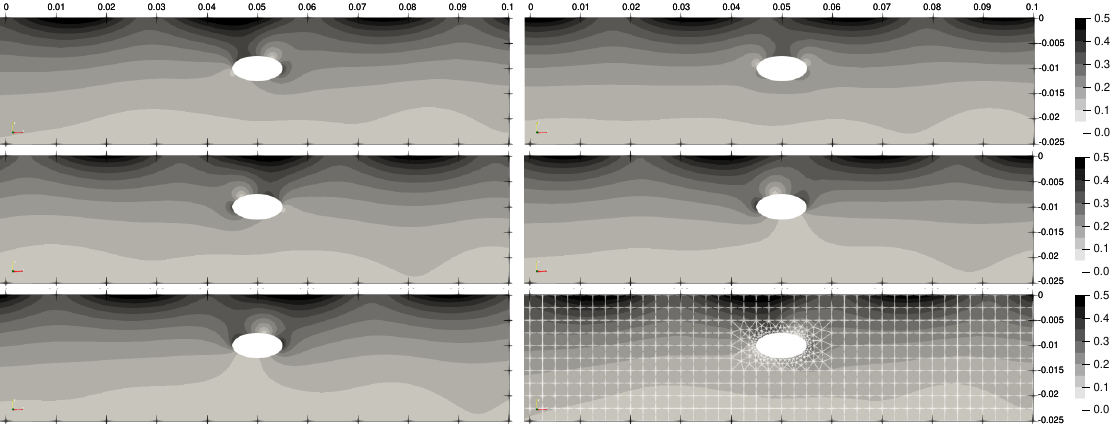}
	\caption{Special test - $||\nabla \varphi||$ for $t=2.20, 2.24, 2.28~s$ (left column from top to bottom) and $t=2.32, 2.36, 2.40~s$ (right column from top to bottom), with the mesh in the last image.}
	\label{fig:videoderniertest}
\end{figure}

The right part of the mesh is shown in Figure \ref{fig:videoderniertest} with the right last image. The rest of the mesh uses the same space discretization.
\begin{figure}[h]
	\centering
	\includegraphics[width=11.9cm]{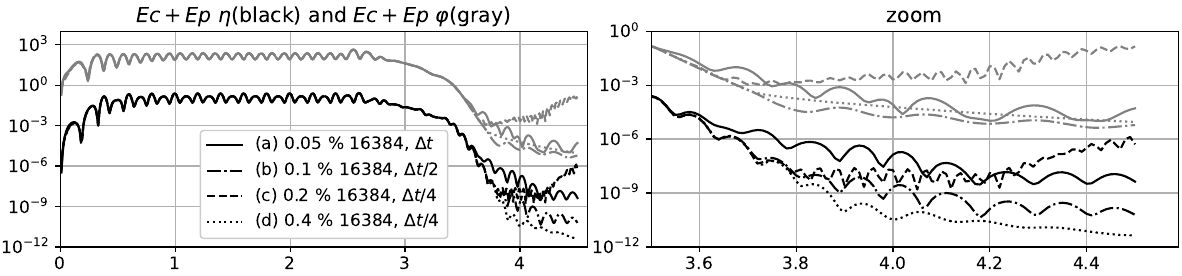}
	\caption{Function of time, the left graph gives the total energy for $\eta$ in $\Gamma_s$ (black color) and for $\varphi$ in $\Omega$ (gray color) for $x > 0$ for the cases (a-b-c-d). The right graph is a zoom of the left graph.}
	\label{fig:videoderniertest_energy}
\end{figure}
Figure \ref{fig:videoderniertest_energy} shows the time evolution of the total energy for just the half right of the domain ($0<x<0.1$) for four different cases. The case (a) uses just $8$ last terms of the Padé approximation ($0.05 \%$ of the ABC order $16384$) with a time step equal to $\Delta t$ (given in the caption of Table \ref{table-cas}). The cases (b), (c) and (d) use respectively $16$ last terms ($0.1 \%$ of $16384$) with a time step equal to $\Delta t/2$, $32$ last terms with $\Delta t/4$ and $64$ last terms with $\Delta t/4$.
For the cases (a), (b) and (d), we keep similar efficient results as previous simulations. The energies decrease to the equivalent levels after an important excitation duration equal to $2.4~s$, between $10^{-6}$ and $10^{-10}$ ($t \simeq 4~s$) for the total energy associated to the surface problem and between $10^{-4}$ and $10^{-5}$ for the basin problem. For the case (a), we use very few terms. With the others cases, we can note that the reduction of the time step must be offset by an increase in the rate of the Padé reduction. But as we can see with the case (c), the ABC does not work well, the energies increase ($t=4.2~s$). The linear correction is not enough. We need to get closer to degree 2 associated to the order two of the time derivative in our model. In other words, to get efficient and low-cost results, it's important to use compatible space and time discretizations associated to (\ref{Pade-reduc-0}). And if we use smaller time steps, the efficient threshold must be change, we must use more last terms for the Padé approximation.\\

For the positive domain ($x \geq 0$), Figure \ref{fig:videoderniertest} shows the magnitude of $\nabla \varphi$ for different times $t \in [2.2, 2.4]$ (the last excitation period).

\section{Conclusion}

For hydrodynamic wave problems with a large number of small objects close to the surface, here just below, the surface tension effect (STE) can be useful. Especially with STE, we have proposed to use a Coupled Hydrodynamic Wave Model here without flow with a small added mass on the surface, this evolution of the Neumann-Kelvin Model allows to use the background of the Absorbing Boundary Conditions community. So, we have developed a specific Coupled Padé-type High-Order Boundary Condition with a compatibility constraint managing the coupling between the ABC associated to the surface model with STE and the ABC of the fluid model below the surface. This compatibility condition uses new coefficients introduced in the Padé-type HABC, similar as Higdon coefficients, avoiding instabilities. These developments are made possible because in this context (presence of multiple objects), FEM is competitive. To explore problems with a big fluid celerity as water, we have proposed to use an important and efficient reduction of the Padé approximation, for that we need to use a sufficient small, not too big added mass on the surface and compatible space, time discretizations or larger time steps (with implicit schemes). We can note that the added mass can be very small, practically negligible in the model. \\

All our numerical computations give controlled accuracy results, with a precision associated to the order of Padé-type HABC and, in particular for water case, a numerical cost associated to the reduction percentage of the number of the terms of the Padé approximation. To allow to get the best chosen precision, the space discretization is important, it's associated to the improvement of the relation with the continuous model by the way of the approximation of the Dirichlet to Neumann operator. \\

Our Padé-type HABC is determined to go out the propagative waves. It's not done for the evanescent waves. It can be a complementary objective, interesting for example to limit a little more the size of the truncated domain. The possible important accuracy of our Padé-type HABC makes it possible to consider its exploitation to construct a very efficient interface condition for Domain Decomposition Method as \cite{Wilk-Jacques-2021}. And the "FEM" coupled model with our CP-HABC may be a way to prospect fluid-structure interactions with multiple objects close to the free surface (\cite{melange-spheres-Valettea}, \cite{microcapsules-DeVuyst}).

\bibliographystyle{plain}

\appendix

\end{document}